\documentclass[11pt]{amsart}

\usepackage{amsmath}
\usepackage{amsthm}
\usepackage{amssymb}
\usepackage{hyperref}
\usepackage[margin=1.3in]{geometry}

\numberwithin{equation}{section}

\newtheorem{prop}{Proposition}
\newtheorem{lemma}[prop]{Lemma}
\newtheorem{thm}[prop]{Theorem}

\newtheorem*{thm*}{Theorem}
\numberwithin{prop}{section}

\theoremstyle{definition}

\newcommand{\dt}{\frac{\partial}{\partial t}}

\renewcommand{\bar}[1]{\overline{#1}}
\newcommand{\del}{\partial}
\newcommand{\delb}{\bar{\partial}}

\newcommand{\bj}{\bar{j}}

\newcommand{\bl}{\bar{l}}

\newcommand{\bp}{\bar{p}}

\newcommand{\til}[1]{\widetilde{#1}}

\DeclareMathOperator{\Ric}{Ric}
\DeclareMathOperator{\Rm}{Rm}

\DeclareMathOperator{\tr}{tr}

\DeclareMathOperator{\Vol}{Vol}

\begin{document}
	\title{The Behavior of Chern Scalar Curvature under Chern-Ricci Flow}
\author{Matthew Gill}
\address{Department of Mathematics, University of California, Berkeley, 970 Evans Hall \#3840, Berkeley, CA 94720-3840 USA}
\thanks{Supported by NSF RTG grant DMS-0838703.}

\author{Daniel Smith}
\address{Department of Mathematics, Furman University, 3300 Poinsett Highway, Greenville, SC, 29613 USA}
\thanks{}

	\maketitle
	
\begin{abstract}
In this note we study finite-time singularities in the Chern-Ricci flow. We show that finite-time singularities are characterized by the blow-up of the scalar curvature of the Chern connection.
\end{abstract}

\section{Introduction}
Let $(M,g,J)$ be a complex manifold such that $J$ is compatible with the Riemannian metric $g$. Associated to $g$, we have the Hermitian 2-form $\omega$ given by $\omega(X,Y) = g(J(X),Y)$ for vector fields $X$ and $Y$. We call a manifold admitting such a structure a \emph{Hermitian manifold}. Recently a number of geometric flows of Hermitian metrics have been introduced to study the underlying Hermitian structure of complex manifolds (e.g. \cite{Gill1}, \cite{CRF1}, \cite{PCF1}, \cite{HCF}, and \cite{Liu-Yang}). The Chern-Ricci flow was introduced by the first named author in \cite{Gill1} and further developed by Tosatti-Weinkove in \cite{CRF1}. The Chern-Ricci flow is written
\begin{gather} \label{CRF}
\begin{split}
	\dt \omega &= -\Ric^C(\omega) \\
	\omega(0) &= \omega_0
\end{split}
\end{gather}
where $\Ric^C(\omega) = -\sqrt{-1} \del \bar{\del} \log \det (g)$. In other words, $\Ric^C(\omega)$ is the curvature of the anti-canonical bundle. Notice that when $\omega$ is a K\"{a}hler metric, $\Ric^C(\omega)$ coincides with the Ricci curvature of $\omega$. Hence if $\omega_0$ is K\"{a}hler, then (\ref{CRF}) is the K\"{a}hler-Ricci flow.\\
\indent The Chern-Ricci flow is a diffusion-reaction equation and therefore one expects the development of singularities in the flow. In fact, explicit examples of singularity development are described by Tosatti-Weinkove in \cite{CRF1}, \cite{CRF2}, and Tosatti-Weinkove-Yang in \cite{CRF3}. Moreover in \cite{CRF1} Tosatti-Weinkove give characterizations of finite-time singularities. To describe these characterizations, first notice that the Chern-Ricci flow (\ref{CRF}) my be written
\begin{align*}
	\dt \omega &= -\Ric^C(\omega_0) + \sqrt{-1}\del \bar{\del} \theta(t) \\
	\omega(0) &= \omega_0
\end{align*}
where $\theta(t) = \log \frac{\det(g(t))}{\det(g_0)}$. Hence the solution $\omega(t)$ can be expressed as
\begin{align}
	\omega(t) = \alpha_t + \sqrt{-1}\del \bar{\del}\phi(t), \text{ where } \alpha_t = \omega_0 - t\Ric^C(\omega_0)
	\label{CRF solution}
\end{align} 
and $\phi(t)$ is a function satisfying $\dt \phi(t) = \theta(t)$ and $\phi\big|_{t=0} = 0$. Let
\begin{align*}
	T = \sup \{ t \geq 0 |\text{ $\exists \psi \in C^{\infty}(M)$ with $\alpha_t + \sqrt{-1}\del \bar{\del} \psi > 0$}\}.
\end{align*}
It is clear that the flow cannot exist beyond time $T$. It was shown by Tosatti-Weinkove that the flow (\ref{CRF}), in fact, exists up to time $T$.
\begin{thm*} (Tosatti-Weinkove \cite{CRF1})
	There exists a unique maximal solution to the Chern-Ricci flow on $[0,T)$.
\end{thm*}
In addition, Tosatti-Weinkove give the following description of finite-time singularities on complex surfaces.
\begin{thm*} (Tosatti-Weinkove \cite{CRF1})
	Let $M$ be a compact complex surface and $\omega_0$ a $\del \bar{\del}$-closed Hermitian metric. Then the Chern-Ricci flow starting at $\omega_0$ exists until either the volume of $M$ goes to zero, or the volume of a curve of negative self-intersection goes to zero.
\end{thm*}
In this note we give an alternative characterization of finite-time singularities in terms of the blow-up of the scalar curvature. In particular we prove
\begin{thm} \label{scalar curvature blow-up}
	Let $M$ be a compact complex manifold of complex dimension $n$ and $\omega_0$ a Hermitian metric.
	Then the solution $\omega(t)$ of the Chern-Ricci flow (\ref{CRF}) exists on the maximal interval $[0,T)$ and 
	either $T = \infty$ or
		$$\limsup_{t \rightarrow T} \left(\sup_{M}R(g(t)) \right) = \infty$$
	where $R$ denotes the scalar curvature of the Chern connection.
\end{thm}
This is a generalization of Zhang's analogous result for K\"{a}hler-Ricci flow \cite{Zhang1} and the proof employs similar arguments. Tosatti-Weinkove-Yang have also investigated the Chern scalar curvature along the normalized Chern-Ricci flow in the case of elliptic surfaces (with genus strictly larger than one). In particular, they show that the Chern scalar curvature is bounded above by $Ce^{t/2}$ and below by $-C$ and apply this estimate to show that the flow converges exponentially fast to the base in the Gromov-Hausdorff sense \cite{CRF2}

Now we specialize to the case when $(M, g_0)$ is a compact Gauduchon surface, i.e. the associated $(1,1)$ form $\omega_0$ of $g_0$ has the property that
\begin{equation}
\del \delb \omega_0 = 0.
\end{equation} 
A theorem of Gauduchon shows that every Hermitian metric on a complex surface is conformal to a unique Gauduchon metric. When the Chern-Ricci flow terminates at a finite time $T$, there are two cases to consider: finite time collapsing and non-collapsing. If $\Vol(M,g(t)) \to 0$ as $t \to T$, we say that the flow is \emph{collapsing}, otherwise we say it is \emph{non-collapsing}. Tosatti-Weinkove consider the case of finite time non-collapsing on compact Gauduchon surfaces and show that there are finitely many $(-1)$ curves on $M$ away from which the flow converges smoothly on compact sets \cite{CRF2}.

Using the proof of Theorem \ref{scalar curvature blow-up} and their result, we are able to show that the set on which singularities develop on $M$ is equal to the region on which Chern scalar curvature blows up.  Collins-Tosatti have an analogous result in the K\"ahler case \cite{CT}. Let $g(t)$ be the solution to the Chern-Ricci flow that is non-collapsing in finite time. Let $E$ denote the union of the $(-1)$ curves on $M$ that are collapsing as $t \to T$. Define the singularity formation set $\Sigma$ to be the compliment of all points $x$ in $M$ such that there exists an open neighborhood $U$ of $x$ on which $|\Rm(t)| \leq C$ for all $0 \leq t < T$. Also define the set $\Sigma'$ to be the compliment of all points $x$ in $M$ such that there exists an open neighborhood $U$ of $x$ on which $|R(t)| \leq C$ for all $0 \leq t < T$. Then we have the result:
\begin{thm} \label{singularity set}
Let $(M,g_0)$ be a compact Gauduchon surface that is non-collapsing in finite time along the Chern-Ricci flow. Then
\begin{equation*}
E = \Sigma = \Sigma'.
\end{equation*}
If instead the flow is collapsing in finite time, we have
\begin{equation*}
M = \Sigma = \Sigma'.
\end{equation*}
\end{thm}

In the final section we include an explicit example of a Type I singularity, in the sense that the Chern scalar curvature blows up like $C(T-t)^{-1}$, along the Chern-Ricci flow in the case of the Hopf manifold.

\section{Preliminaries}
Let $(M,g,J)$ be a closed Hermitian manifold of complex dimension $n$. The \emph{Chern connection} is the unique connection $\nabla$ satisfying: $\nabla \omega = 0$, $\nabla J = 0$, and $T^{1,1} = 0.$ In local holomorphic coordinates, $z_1, \dots, z_n$, the only non-vanishing components of the Chern connection are
\begin{align}
	\Gamma_{ij}^k = g^{k \bl} \del_i g_{j \bl} 
\end{align}
and its conjugate. The curvature of the Chern connection is given locally by
\begin{align}
	R_{i \bj k}^{l} = -\del_{\bj} \Gamma_{ik}^l.
\end{align}
Taking a trace of the Chern curvature over the last two components yields the Chern-Ricci curvature, that is
\begin{align}
	R_{i \bj}^C = g^{k\bl}R_{i \bj k \bl} = - \del_{i} \del_{\bj} \log \det (g).
\end{align}
Interestingly, the parabolic flow of Hermitian metrics defined by Streets-Tian in \cite{HCF} is the trace of the Chern curvature over the first two components plus a term that is quadratic in the torsion of $\nabla$. Finally, the \emph{Chern scalar curvature}, denoted $R$, is defined as the trace of the Chern-Ricci curvature
	$$R = g^{i \bj}R_{i\bj}^C.$$
\indent One of the reasons for the recent success of the Chern-Ricci flow is that, like K\"{a}hler-Ricci flow, the Chern-Ricci flow is equivalent to a scalar flow. To see this, first recall from (\ref{CRF solution}) that the solution of the Chern-Ricci flow can be written
\begin{align}
	\omega(t) = \alpha_t + \sqrt{-1}\del \bar{\del} \phi(t). \label{CRF solution 2}
\end{align}
Differentiating (\ref{CRF solution 2}) in time, we have
\begin{align}
	-\Ric^C(t) = \sqrt{-1}\del \bar{\del} \log \det(g_0) + \sqrt{-1}\del \bar{\del}\left( \dt \phi(t) \right). \label{equation}
\end{align}
Rearranging terms and using (\ref{CRF solution 2}) yields
\begin{align*}
	\sqrt{-1}\del \bar{\del}\left( \dt \phi(t) \right) &= \sqrt{-1}\del \bar{\del} \left( \log \det(g(t)) - \log \det(g_0) \right) \\
		&= \sqrt{-1}\del \bar{\del} \log  \frac{(\alpha_t + \sqrt{-1}\del \bar{\del}\phi(t))^n}{\omega_0^n}.
\end{align*}
Therefore the Chern-Ricci flow is equivalent to the scalar complex Monge-Amp\`{e}re equation
\begin{align}
	\dt \phi &= \log \frac{(\alpha_t + \sqrt{-1}\del \bar{\del}\phi(t))^n}{\omega_0^n}, \text{ with } \phi\big|_{t=0} = 0. \label{CX MA}
\end{align}
\indent The proof of Theorem \ref{scalar curvature blow-up} will repeatedly exploit the fact that the Chern-Ricci flow (\ref{CRF}) is equivalent to the scalar flow (\ref{CX MA}).

\section{Proof of the Theorem \ref{scalar curvature blow-up}}
In this section we show that finite-time singularities of the Chern-Ricci flow occur if and only if the scalar curvature is unbounded as $t$ approaches the singular time. 
\proof We will proceed in the same fashion as Zhang in \cite{Zhang1} (cf. Section 7.2 of \cite{Song-Weinkove1}). Let $\omega(t)$ be a solution of the Chern-Ricci flow (\ref{CRF}) on the maximal interval $[0,T)$ with $T < \infty$. Suppose, by way of contradiction, that there exists a constant $C$ such that
\begin{align}
	\limsup_{t \rightarrow T} \left(\sup_{M}R(g(t)) \right) \leq C. \label{R upper bound}
\end{align}
Throughout this section $C$ will denote a constant which is independent of time, however this constant may change from line to line. We will show that if (\ref{R upper bound}) holds, then for any $k \geq 0$ we can produce uniform $C^k$ bounds on the metric $\omega(t)$ and hence $T$ is not maximal. We begin by obtaining uniform $C^0$ estimates on $\omega(t)$ for $t \in [0,T)$.\\
\indent First, by studying the evolution of the scalar curvature $R(t)$ under (\ref{CRF}), we will show that the $C^0$ norm of $R(t)$ is uniformly bounded on $[0,T)$.
\begin{lemma} \label{evolution of R}
Under the Chern-Ricci flow (\ref{CRF}) the scalar curvature of the Chern connection evolves by
$$\dt R = \Delta R + |\Ric^C(\omega)|^2.$$
\end{lemma}
\proof In local holomorphic coordinates, at a point $p \in M$, we compute
\begin{align*}
	\dt R &= \dt \left( g^{i\bj}R_{i\bj}^C \right) \\
		&= g^{k\bj}g^{i\bp}R_{k\bp}^CR_{i\bj}^C - g^{i\bj}\del_{i}\del_{\bj}\left( g^{k\bl}\dt g_{k\bl} \right)  \\
		&= \Delta R + |\Ric^C(\omega)|^2.
\end{align*}
\qed
\medskip

Hence by Lemma \ref{evolution of R} and the maximum principle, it follows that $R(t)$ is uniformly bounded below on $[0,T)$. Combining this lower bound with (\ref{R upper bound}), we have that $||R(t)||_{C^0}$ is uniformly bounded on $[0,T)$.\\
\indent Next we will use the bounds on the scalar curvature $R(t)$ to obtain bounds on the determinant of $\omega(t)$, which are independent of $t$. Recall from (\ref{CX MA}) that the Chern-Ricci flow is equivalent to the complex Monge-Amp\`{e}re equation $\dt \phi = \log \left( \frac{\omega^n}{\omega_0^n} \right)$. Now, since $\dt \log (\omega^n) = -R$, we have
\begin{align}
	\dt \dot{\phi} = -R. \label{evolution of phi dot}
\end{align}
Consequently,
\begin{align}
	\left|\dt \dot{\phi}\right| = |R| \leq C. \label{phi dot upper bound}
\end{align}
Integrating (\ref{phi dot upper bound}) in time we have that there exists a constant $C$ such that 
\begin{align}
	\left|\log\left( \frac{\omega^n}{\omega_0^n} \right)\right| = \left|\dt\phi\right| \leq C. \label{determinant bound}
\end{align}
Notice that (\ref{determinant bound}) gives uniform bounds on the determinant of $\omega(t)$ on $[0,T)$. Integrating in time again gives uniform bounds for $|\varphi|$. Given uniform estimates on the determinant of the metric, to prove that $\omega(t)$ is uniformly bounded in $C^0$, it suffices to show that the trace of the metric $\tr_{g_0}g$ is bounded above.\\
\indent For a bound on the trace of the metric, we repeat the estimates derived by Tosatti-Weinkove in Section 4 of \cite{CRF1} with minor modifications. First let 
	$$Q_1 \doteq t\dot{\phi} - \phi -nt$$ 
and notice that $Q_1$ is uniformly bounded. Combining the trace of (\ref{equation}) with (\ref{evolution of phi dot}), it follows that
\begin{align*}
	\left( \dt - \Delta \right) \dot{\phi} = \tr_{\omega}\left(-\Ric^C(\omega_0)\right).
\end{align*}
Taking the trace of (\ref{CRF solution 2}), we have $\Delta \phi = \tr_{\omega}(t\Ric^C(\omega_0) - \omega_0) + n$ and so,
\begin{align}
	\left( \dt - \Delta \right)Q_1 = -\tr_{\omega}\omega_0. \label{heat operator Q1}
\end{align}
Let $\til{C}$ denote a constant which is sufficiently large so that $\phi + \til{C} \geq 1$ and let $B$ be a large constant to be chosen below. Next, define
\begin{align}
	Q_2 \doteq \log \tr_{g_0}g - \phi + \frac{1}{\phi + \til{C}} + BQ_1. \label{Q2}
\end{align}
Notice that each term in $Q_2$, with the possible exception of $\log \tr_{g_0}g$, is uniformly bounded. Therefore it is enough to show that $Q_2$ is bounded above at a point where $Q_2$ achieves its maximum.\\
\indent By Proposition 3.1 in \cite{CRF1}, we have the following estimate on the heat operator applied to $\log \tr_{g_0}g$.
\begin{prop} (Tosatti-Weinkove \cite{CRF1})
	Given a solution $g(t)$ of the Chern-Ricci flow (\ref{CRF}),
	\begin{align}
		\left( \dt - \Delta \right) \log \tr_{g_0}g \leq \frac{2}{(\tr_{g_0}g)^2}Re \left( g^{\bl k}(T_0)_{kp}^p
		\del_{\bl}\tr_{g_0}g \right) + C\tr_{g}g_0. \label{log trace bound}
	\end{align}
\end{prop}
Now, we may assume that at a maximum of $Q_2$ that $\tr_{g_0}g \geq 1$ since otherwise we have the desired upper bound. For a bound on the first term in (\ref{log trace bound}), notice that at a maximum of $Q_2$, $\del_i Q_2 = 0$. Equivalently,
\begin{align}
	\frac{1}{\tr_{g_0}g}\del_i \tr_{g_0}g - \del_i \phi - \frac{1}{(\phi + \til{C})^2}\del_i \phi = 0. \label{critical point}
\end{align}
Hence at a maximum of $Q_2$,
\begin{align}
	& \left| \frac{2}{(\tr_{g_0}g)^2}Re \left( g^{\bl k}(T_0)_{kp}^p \del_{\bl} \tr_{g_0}g \right) \right| \nonumber \\
		&\leq \left| \frac{2}{\tr_{g_0}g}Re \left( \left( 1 + \frac{1}{(\phi + \til{C})^2} \right) g^{\bl k}(T_0)_{kp}^p 
		(\del_{\bl} \phi) \right) \right| \nonumber \\
		&\leq \frac{|\del \phi|_{g}^2}{(\phi + \til{C})^3} + C(\phi + \til{C})^3 \frac{\tr_{g}g_0}{(\tr_{g_0}g)^2}. \label{trace bds}
\end{align}
Notice that we may also assume that at a maximum of $Q_2$, $(\tr_{g_0}g)^2 \geq (\phi + \til{C})^3$ since if this is not the case then we again have the desired upper bound. And so, by (\ref{heat operator Q1}), (\ref{log trace bound}), and (\ref{trace bds}), at a maximum
\begin{align*}
	0 	&\leq \left( \dt - \Delta \right) Q_2 \\
		&\leq \frac{|\del \phi|_{g}^2}{(\phi + \til{C})^3} + C \tr_{g}g_0 - \left( 1 + \frac{1}{(\phi + \til{C})^2} \right) \dot{\phi}
			+ \left( 1 + \frac{1}{(\phi + \til{C})^2} \right) \tr_{g}(g - \alpha_t) \\
		&\hspace{.18 in} - \frac{2}{(\phi + \til{C})^3} |\del \phi|_{g}^2 - B \tr_{g}g_0 \\
		&\leq (C - B) \tr_{g}g_0 + C'
\end{align*}	
where $C'$ comes from the bounds on $\phi$ and $\dot{\phi}$. Also notice that we used the bound
\begin{align*} 
	-\tr_{g}\alpha_t &= -\tr_{g}g_0 + t\tr_g\left( \Ric^C(g_0) \right) \\
		&\leq -\tr_{g}g_0 + \hat{C}\tr_{g}g_0 
\end{align*}
where $\hat{C}$ is a constant that depends on the geometry of $g_0$ and harmlessly on $t$. Now we choose $B = C + 1$. Hence,
\begin{align}
	\tr_{g}g_0 \leq C'. \label{wrong trace bound}
\end{align}
Finally, using the bounds on the determinant of the metric (\ref{determinant bound}) and (\ref{wrong trace bound}),
\begin{align}
	\tr_{g_0}g \leq \frac{1}{(n - 1)!} (\tr_{g}g_0)^{n-1} \frac{\det g}{\det g_0} \leq C''.
\end{align}
Thus we have shown that $Q_2$ is uniformly bounded above and so $\tr_{g_0}g$ is bounded above. This proves the uniform $C^0$ bounds on $\omega(t)$; hence there exists a constant $C$ so that
\begin{align*}
	\frac{1}{C}\omega_0 \leq \omega(t) \leq C \omega_0
\end{align*}
for $t \in [0,T)$. \\
\indent The higher order estimates on the metric follow from the first named author in \cite{Gill1}. This proves the theorem.
\qed 

\section{Proof of the Theorem \ref{singularity set}}

The proof follows rather quickly from the previous proof and the results of Tosatti-Weinkove and follows a method similar to that in \cite{CT}. In the case of finite time non-collapsing, Tosatti-Weinkove have shown that there exist finitely many disjoint $(-1)$ curves $E_1, \ldots, E_l$ on $M$ with $\Vol(E_i, g(t)) \to 0$ as $t \to T$ \cite{CRF1}. This immediately implies that $E \subset \Sigma$. Tosatti-Weinkove also proved that the flow converges smoothly on compact sets away from $E$, giving the reverse containment since no singularities develop outside of $E$ and hence $E = \Sigma$ \cite{CRF2}. 

Trivially, we have $\Sigma' \subset \Sigma = E$. It remains to show that for any $x$ in $M$ with an open neighborhood $U$ on which $|R(t)| \leq C$ for all $0 \leq t < T$, then we have smooth convergence on a compact set $K \subset U$. But this is exactly what was proved in the previous section if we instead perform all calculations on $K$ instead of $M$. The case of finite time collapsing follows analogously. 
\qed

\section{The Hopf Manifold}

Let $\alpha = (\alpha_1, \ldots, \alpha_n) \in \mathbb{C}^n \setminus \{0\}$ with $|\alpha_1| = \ldots = |\alpha_n| \neq 1$. Then consider the Hopf manifold $M_\alpha = ( \mathbb{C} \setminus \{0\} ) / \sim$ where $(z_1, \ldots, z_n) \sim (\alpha_1 z_1, \ldots, \alpha_n z_n)$. Let
\begin{align}
& \omega_H = \frac{\delta_{i\bj}}{r^2} \sqrt{-1} dz^i \wedge dz^{\bj} \\
& \Ric^C(\omega_H) = \frac{n}{r^2} \left( \delta_{i\bj} - \frac{z_i z_{\bj}}{r^2} \right) \sqrt{-1} dz^i \wedge dz^{\bj}.
\end{align}
Tosatti-Weinkove show that the metric
\begin{equation}
\omega(t) = \omega_H - t \Ric^C(\omega_H) = \left( \frac{\delta_{i{\bj}}(1-nt)}{r^2} - \frac{nt z_i z_{\bj}}{r^4} \right) \sqrt{-1} dz^i \wedge dz^{\bj}
\end{equation}
is a solution to the Chern-Ricci flow on the Hopf manifold \cite{CRF1}. In particular, they show that $\Ric^C(\omega(t)) = \Ric^C(\omega_H)$. 

Computing directly, we see that
\begin{equation}
R(\omega(t)) = \tr_{\omega(t)} \Ric^C(\omega_H) \leq \frac{C}{\frac{1}{n} - t}
\end{equation}
and hence the flow develops a Type I singularity at $T = \frac{1}{n}$.

\section{Acknowledgments}

The authors would like to especially thank Valentino Tosatti for pointing out that Theorem \ref{singularity set} should follow from the proof of Theorem \ref{scalar curvature blow-up}. The second named author would also like to thank Jon Wolfson for numerous helpful conversations.

\bibliographystyle{plain}
\bibliography{biblio}

\begin{thebibliography}{10}

\bibitem{CT}
T.~{Collins} and V.~{Tosatti}.
\newblock {K\"{a}hler currents and null loci}.
\newblock {\em ArXiv e-prints}, May 2013.

\bibitem{Gill1}
M.~{Gill}.
\newblock Convergence of the parabolic complex {M}onge-{A}mp\`ere equation on
  compact {H}ermitian manifolds.
\newblock {\em Comm. Anal. Geom.}, 19(2):277--303, 2011.

\bibitem{Liu-Yang}
K.~{Liu} and X.~{Yang}.
\newblock {Geometry of Hermitian manifolds}.
\newblock {\em ArXiv e-prints}, October 2010.

\bibitem{Song-Weinkove1}
J.~{Song} and B.~{Weinkove}.
\newblock {Lecture notes on the K\"{a}hler-Ricci flow}.
\newblock {\em ArXiv e-prints}, December 2012.

\bibitem{PCF1}
J.~{Streets} and G.~{Tian}.
\newblock A parabolic flow of pluriclosed metrics.
\newblock {\em Int. Math. Res. Not. IMRN}, (16):3101--3133, 2010.

\bibitem{HCF}
J.~{Streets} and G.~{Tian}.
\newblock Hermitian curvature flow.
\newblock {\em J. Eur. Math. Soc. (JEMS)}, 13(3):601--634, 2011.

\bibitem{CRF1}
V.~{Tosatti} and B.~{Weinkove}.
\newblock {On the evolution of a Hermitian metric by its Chern-Ricci form}.
\newblock {\em ArXiv e-prints}, December 2012.

\bibitem{CRF2}
V.~{Tosatti} and B.~{Weinkove}.
\newblock {The Chern-Ricci flow on complex surfaces}.
\newblock {\em ArXiv e-prints}, September 2012.

\bibitem{CRF3}
V.~{Tosatti}, B.~{Weinkove}, and X.~{Yang}.
\newblock {Collapsing of the Chern-Ricci flow on elliptic surfaces}.
\newblock {\em ArXiv e-prints}, February 2013.

\bibitem{Zhang1}
Z.~{Zhang}.
\newblock Scalar curvature behavior for finite-time singularity of
  {K}\"ahler-{R}icci flow.
\newblock {\em Michigan Math. J.}, 59(2):419--433, 2010.

\end{thebibliography}

\end{document}